\documentclass[11pt]{article}
\usepackage{amsmath, amsfonts, amssymb}
\usepackage{graphicx, amsthm,color}
\usepackage{cite}
\usepackage[english]{babel}
\usepackage[utf8]{inputenc}
\usepackage{amscd}
\usepackage{newlfont}
\usepackage{enumerate}

\newtheorem{theorem}{Theorem}[section]
\newtheorem{lemma}[theorem]{Lemma}

\newtheorem{proposition}[theorem]{Proposition}
\theoremstyle{definition}

\newtheorem{example}[theorem]{Example}
\newtheorem{definition}[theorem]{Definition}
\theoremstyle{remark}
\newtheorem{remark}[theorem]{Remark}
\numberwithin{equation}{section}

\begin{document}

\begin{center}
{\LARGE Extending Rademacher Theorem to Set-Valued Maps.}
\end{center}

\medskip

\begin{center}
{\large \textsc{Aris Daniilidis \& Marc Quincampoix}}
\end{center}

\bigskip

\noindent\textbf{Abstract.} Rademacher theorem asserts that Lipschitz
continuous functions between Euclidean spaces are differentiable almost
everywhere. In this work we extend this result to set-valued maps using an
adequate notion of set-valued differentiability relating to convex processes.
Our approach uses Rademacher theorem but also recovers it
as a special case. \bigskip


\vspace{0.6cm}

\textit{Secondary}: 49J53, 49J52

\section{Introduction}

Rademacher theorem asserts that a (locally) Lipschitz continuous
function $f:\mathbb{R}^{d}\rightarrow\mathbb{R}^{\ell}$ is differentiable
almost everywhere (that is, the derivative $df(x)$ of $f$ exists at every
point $x\in\mathbb{R}^{d}\setminus\mathcal{N}$ where $\mathcal{N}$ has
Lebesgue measure zero). In this paper we study extensions of this result to
Lipschitz continuous set-valued maps $F$ from $\mathbb{R}^{d}$ to
$\mathbb{R}^{\ell}$ with nonempty convex compact values. We recall
(see~\cite{AF1990-book}, \cite{Morduk2006-book}, \cite{RW1998-book}
\emph{e.g.}) that a set-valued map $F$ is called locally Lipschitz at $\bar
{x}$ if for some $k\geq0$ and all $x,x^{\prime}$ in a neighborhood of $\bar
{x}$ it holds:
\[
F(x)\subset F(x^{\prime})+k\,|x-x^{\prime}|\,\mathrm{B},
\]
where $\mathrm{B}$ is the closed unit ball of $\mathbb{R}^{\ell}$.  A special type of set-valued maps
are the so-called \textit{convex processes} (see \cite{R1972}, \cite{R1967}). A set-valued map $T:\mathbb{R}^{d}
\rightrightarrows\mathbb{R}^{\ell}$ is called a convex process, if its graph is a closed convex cone in
$\mathbb{R}^{d}\times\mathbb{R}^{\ell}$, or equivalently, if $0\in T(0)$ and for every
$x,y\in\mathbb{R}^{d}$ and $\lambda>0$ it holds:
\begin{equation}
T(x+y)\subset T(x)+T(y)\qquad\text{and }\qquad T(\lambda x)=\lambda T(x).
\label{eq:convex-process}
\end{equation}
Convex processes are relevant in many applications in control and optimization
(see \cite{ACG2005}, \cite{AF1986}, \cite{L1994}~\emph{e.g.} or the classical monographs
\cite{AF1990-book}, \cite{Morduk2006-book}, \cite{RW1998-book}). They are considered as set-valued analogues of linear continuous operators; we
refer the reader to \cite[Chapter 2]{AF1990-book} for a discussion on this fact. We hereby use this notion to define
differentiability of a set-valued map: indeed, we say that $F$ is differentiable at $(x,y)\in\mathbb{R}^{d}\times
\mathbb{R}^{\ell}$ with $y\in F(x)$ if the graphical derivative of $F$ there
is a \textit{convex process} (see forthcoming Definition~\ref{diffF}).\smallskip\newline
Let us mention for
completeness that an alternative notion of differentiability of set-valued
maps, called $H$-differentiability, has been introduced in \cite{P2011} and
studied in \cite{DP2011}, \cite{GGJ2015}. The notion of $H$-differentiability
is based on the (less restrictive) positively homogeneous set-valued operators
(which are not necessarily convex processes) and its main drawback is that
leads to conic approximations of $F$ that need not be unique (see discussion
in \cite[Section 2.3]{DP2011}). \medskip\newline
The manuscript is organized as follows: in the next section we fix our
notations and give the main definitions and preliminary results.
Section~\ref{sec:3} is dedicated to differentiability of set-valued maps and
its relation with the case of differentiability of functions. In
Section~\ref{sec:4} we obtain an extension of Rademacher theorem for the class
of set-valued maps that are generated by a finite number of pointwise affinely
independent Lipschitz functions (Proposition~\ref{cor:family}). This result recovers the classical Rademacher
theorem as a special case. In Section~\ref{sec:5} we obtain a general
differentiability result (Theorem~\ref{RadeFint}) for the class of isotropically Lipschitz functions (Definition~\ref{def_isotropic}) with convex compact values with nonempty
interior and smooth boundary.

\section{Preliminaries}

\label{sec:2} Throughout this work, we denote by $|x|$ the Euclidean norm of
an element $x\in\mathbb{R}^{d}$, and by $\mathrm{B}_{X}:=\{x\in\mathbb{R}^{d}:\,|x|\leq1\}$
the closed unit ball, centered at the origin of the normed
space $X$. The index will often be omitted if there is no ambiguity about the
space. We set $B_{\delta}(x):=x+\delta\mathrm{B}_{X}$ for the closed ball
centered at $x$ with radius $\delta>0.$ Given a sequence $(t_{n})_{n\geq1}$ of
real numbers, we shall use the notation $t_{n}\searrow0^{+}$ to indicate that
$t_{n}>0$ and $\underset{n\rightarrow\infty}{\lim}t_{n}=0.$\smallskip

Let further denote by $\mathcal{L}_{d}$ the Lebesgue measure of $\mathbb{R}
^{d}$. A set $A\subset\mathbb{R}^{d}$ is said to have a full measure in
$\mathbb{R}^{d}$ if its complement is contained in a null set, that is,
$\mathbb{R}^{d}\diagdown A\subset\mathcal{N}$ with $\mathcal{L}_{d}
(\mathcal{N})=0$. If a property holds for all points of a full-measure set
$A$, then we say that the property holds \textit{almost everywhere} (in short,
\emph{a.e.}) or \textit{for almost all} $x\in\mathbb{R}^{d}$ (in short,
$\forall_{\mathit{a.e.}}x\in\mathbb{R}^{d}$) and omit the explicit reference
to the set $A$.\smallskip

A function $f:\mathbb{R}^{d}\longrightarrow\mathbb{R}^{\ell}$ is called
locally Lipschitz continuous at $\bar{x}\in\mathbb{R}^{d}$ if there exists
$k\geq0$ such that for all $x,x^{\prime}$ in some neighborhood of $\bar{x}$ it
holds:
\begin{equation}
f(x)-f(x^{\prime})\leq k\,|x-x^{\prime}|. \label{eq:lip}
\end{equation}
We say that $f$ is locally Lipschitz, if it is locally Lipschitz at every
point of its domain. In particular, if \eqref{eq:lip} holds for all
$x,x^{\prime}\in\mathbb{R}^{d},$ then we say that $f$ is $k$-Lipschitz
continuous. We use the notation $df(x)$ to denote the (Fr\'echet) derivative of
$f$ at $x$ (whenever it exists). Notice that $df(x)$ is a linear function from
$\mathbb{R}^{d}$ to $\mathbb{R}^{\ell}$ (linear approximation of $f$ aroung
$x$) and its graph is an affine subspace of $\mathbb{R}^{d}\times
\mathbb{R}^{\ell}$ passing through $(x,f(x)).$\smallskip\newline In this
classical setting, let us recall the following well-known result (see
\cite{Federer1969-book}, \cite{EG1992-book} \emph{e.g.}).\medskip

(\textit{Rademacher theorem}) Every locally Lipschitz function
$f:\mathbb{R}^{d}\rightarrow\mathbb{R}^{\ell}$ is \emph{a.e.} differentiable.

\medskip

\noindent Set-valued maps are tightly related to nonsmooth phenomena and have
been widely used in problems in calculus of variations, control and
optimization (see \cite{GGJ2015}, \cite{GGM2016}, \cite{GS2020}, \cite{Q1992},
\cite{U2021} \emph{e.g.}). For a set-valued map $F$ from $\mathbb{R}^{d}$ to
$\mathbb{R}^{\ell}$, we will use the notation $F:\mathbb{R}^{d}
\rightrightarrows\mathbb{R}^{\ell}$. The graph of $F$ is defined as follows:
\[
\mathrm{gph}(F):=\{(x,y)\in\mathbb{R}^{d}\times\mathbb{R}^{\ell}:\;y\in
F(x)\,\}.
\]

\noindent A set-valued map $F$ is called locally Lipschitz at $\bar{x}$ (or
$k$-Lipschitz around $\bar{x}$) if for some $k\geq0$ and all $x,x^{\prime}$ in
a neighborhood of $\bar{x}$ it holds:
\begin{equation}
F(x)\subset F(x^{\prime})+k\,|x-x^{\prime}|\,\mathrm{B}. \label{eq:Lip}
\end{equation}
Similarly to the case of (single-valued) functions, a set-valued map $F$ is
called locally Lipschitz, if it is locally Lipschitz at every $\bar{x}
\in\mathbb{R}^{d}.$ Moreover, if \eqref{eq:Lip} holds for all $x,x^{\prime}
\in\mathbb{R}^{d},$ then we say that $F$ is $k$-Lipschitz. It is easily seen
that every locally Lipschitz function has closed graph.\smallskip

\noindent Before we proceed further, let us recall that for a closed set
$K\subset\mathbb{R}^{m}$ and $\bar{x}\in K$ the tangent (or contingent) cone
is defined as follows (see \cite{AF1990-book}, \cite{Morduk2006-book},
\cite{RW1998-book} \emph{e.g.})
\begin{equation}
T_{K}(\bar{x}):=\{u\in\mathbb{R}^{m}:\,\exists(h_{n})\searrow0^{+}
,\,\exists(u_{n})_{n}\rightarrow u\,\mbox{ such that }\bar{x}+h_{n}u_{n}\in
K,\,\forall n\geq0\,\}. \label{eq:tangent}
\end{equation}
It follows easily from the definition that $T_{K}(\bar{x})$ is a nonempty
closed cone of $\mathbb{R}^{k}$.\smallskip\newline
 We finally recall that vectors $\{V_{1},\ldots,V_{N}\}\subset
\mathbb{R}^{d}$ are called affinely independent, if for every $\{\mu
_{i}\}_{i=1}^{N}\subset\mathbb{R}$ it holds:
\[
{\displaystyle\sum\limits_{i=1}^{N}} \mu_{i}=0\quad\text{and\quad}
{\displaystyle\sum\limits_{i=1}^{N}} \mu_{i}V_{i}=0\qquad\Longrightarrow
\qquad\mu_{1}=\ldots=\mu_{N}=0.
\]
Notice that this is equivalent to the fact that the family $\{V_{i}-V_{N}\}_{i=1}^{N-1}$ is linearly independent
and yields that $N\leq d+1$.

\section{Graphical derivative and differentiability of set-valued maps}

\label{sec:3}It is natural to apply the notion of tangent cone to the graph of
a set-valued map in order to obtain a linearization of a set-valued map around
a point of its graph. This is achieved by setting $K=\mathrm{gph}(F)$
(assuming closed) and taking the tangent cone at a point $(x,y)\in
\mathrm{gph}(F)$. We formalize this in the following definition (see
\cite{AF1990-book}).

\begin{definition}
[graphical derivative]\label{def_gr-der} Let $F:\mathbb{R}^{d}
\rightrightarrows\mathbb{R}^{\ell}$ be a set-valued map with closed graph. The
graphical derivative of $F$ at $(\bar{x},\bar{y})\in\mathrm{gph}(F)$ is a
set-valued map $DF(\bar{x},\bar{y}):\mathbb{R}^{d}\rightrightarrows
\mathbb{R}^{\ell}$ whose graph is the tangent cone at $(\bar{x},\bar{y})$ of
the graph of $F$, that is:
\[
\mathrm{gph}(DF(\bar{x},\bar{y})):=T_{\mathrm{gph(}F)}(\bar{x},\bar{y}).
\]
\end{definition}

\noindent In view of the above definition, taking into account~\eqref{eq:tangent} we deduce easily that for every
$u\in\mathbb{R}^{d}$:
\begin{equation}
DF(\bar{x},\bar{y})(u):=\{v\in\mathbb{R}^{\ell}:\;\exists(h_{n})\searrow 0^{+},\,\exists(u_{n},v_{n})_{n}\rightarrow(u,v):\;\bar{y}+h_{n}v_{n}\in
F(\bar{x}+h_{n}u_{n}),\,\forall n\,\} \label{DF}
\end{equation}
In case of a locally Lipschitz (set-valued) map, the graphical derivative has
two interesting properties:

\begin{lemma} [graphical derivative of a Lipschitz map] Assume that $F:\mathbb{R}^{d}\rightrightarrows\mathbb{R}^{\ell}$ is $k$-Lipschitz at $\bar{x}$ and pick
any $\bar{y}\in F(\bar{x})$.\smallskip\newline$\mathrm{(i)}$. The graphical derivative
$DF(\bar{x},\bar{y})$ admits a simplified formula:
\begin{equation}
DF(\bar{x},\bar{y})(u):=\{v\in\mathbb{R}^{\ell}:\;\exists(h_{n})\searrow
0^{+},\,\exists(w_{n})_{n}\rightarrow w:\;\bar{y}+h_{n}w_{n}\in F(\bar
{x}+h_{n}u),\,\forall n\,\}. \label{DFk}
\end{equation}
$\mathrm{(ii)}$. For any $u\in\mathbb{R}^{d}$, we have
\[
DF(\bar{x},\bar{y})(u)\,\cap\,k\,|u|\,\mathrm{B}\,\neq\,\emptyset.
\]

\end{lemma}

\noindent\textbf{Proof.} Let us prove (i). A mere comparison with~\eqref{DF}
shows that the right-hand side of~\eqref{DFk} is contained in $DF(\bar{x}
,\bar{y})(u)$. To prove the converse inclusion take $(u,v)\in\mathrm{gph}
\left(  DF(\bar{x},\bar{y})\right)  .$ Then, according to~\eqref{DF}, there
exist sequences $\;\exists(h_{n})\searrow0^{+},$ $(u_{n})_{n}\rightarrow u$
and $(v_{n})_{n}\rightarrow v$ such that $\bar{y}+h_{n}v_{n}\in F(\bar
{x}+h_{n}u_{n})$ for all $n\geq1$. Since $F$ is $k$-Lipschitz, we deduce
from~\eqref{eq:Lip} that $F(\bar{x}+h_{n}u_{n})\subset F(\bar{x}+h_{n}u)+kh_{n}|u-u_{n}|\mathrm{B}$. In particular, there exists a sequence
$(b_{n})_{n}\subset\mathrm{B}$ such that for all $n\geq1$ we have:
\[
\bar{y}+h_{n}(v_{n}+k|u_{n}-u|b_{n})\in F(\bar{x}+h_{n}u).
\]
Set $w_{n}:=v_{n}+k|u_{n}-u|b_{n}$ so that $\bar{y}+h_{n}v_{n}^{\prime}\in
F(\bar{x}+h_{n}u)$ for all $n\geq1.$ Noticing that $(w_{n})_{n}$ converges to
$v$ as $n\rightarrow\infty$ yields \eqref{DFk}.\medskip\newline Let us now
prove (ii). Fix $u\in\mathbb{R}^{d}$ and $(h_{n})\searrow0^{+}$. By Lipschitz
continuity of $F$ we deduce:
\[
\bar{y}\in F(\bar{x})\subset F(\bar{x}+h_{n}u)+\,k\,h_{n}\,|u|\,\mathrm{B}.
\]
Therefore, there exists a sequence $(b_{n})_{n}\subset\mathrm{B}$ such that
$\bar{y}+\,k\,h_{n}\,|u|\,\mathrm{B}\in F(\bar{x}+h_{n}u)$ for all $n\geq1$.
Since $\mathrm{B}$ is compact, passing eventually to a subsequence, we may
assume that $b_{n}\rightarrow b\in\mathrm{B}$ and $w_{n}:=k\,|u|\,b_{n}\rightarrow k\,|u|\,b:=v$. Since $\bar{y}+h_{n}w_{n}\in F(\bar{x}+h_{n}u)$ we
deduce directly from~\eqref{DFk} that $v\in k\,|u|\,\mathrm{B}\cap DF(\bar
{x},\bar{y})(u)$.\medskip\newline The proof is complete. \hfill$\Box$

\bigskip

\noindent Now we are ready to define an appropriate notion of
differentiability of a set-valued map.

\begin{definition}
[Differentiability of a set-valued map]\label{diffF}A set-valued map
$F:\mathbb{R}^{d}\rightrightarrows\mathbb{R}^{\ell}$ with closed graph is said
to be differentiable at $(\bar{x},\bar{y})\in\mathrm{gph}F$ if its graphical
derivative $DF(\bar{x},\bar{y})$ is a convex process (see
\eqref{eq:convex-process}).
\end{definition}

\noindent In other words, $F$ is differentable at $(\bar{x},\bar{y})\in\mathrm{gph}F$ if the closed cone $T_{\mathrm{gph(}F)}(\bar{x},\bar{y})$
is also convex, which means that $DF(\bar{x},\bar{y})$ is a set-valued whose
graph is a closed convex cone.\smallskip\newline Definition~\ref{diffF}
appears here for the first time. The definition is purely geometrical and
guarantees that the derivative is unique (whenever it exists). Moreover, it
generalizes the classical differentiability of Lipschitz functions, as
forthcoming Proposition~\ref{df} will underline. Let us first state the
following lemma.

\begin{lemma}
\label{Lemma-a}Let $f:\mathbb{R}^{d}\longrightarrow\mathbb{R}^{\ell}$ be a
continuous function and assume that for some $\bar{x}\in\mathbb{R}^{d}$ we
have:
\begin{equation}
\displaystyle k_{\bar{x}}:=\limsup_{x\rightarrow\bar{x}}\,\frac{|f(x)-f(\bar
{x})|}{|x-\bar{x}|}<+\infty\label{eq:calm}
\end{equation}
Then setting $F(x)=\{f(x)\}$, for all $x\in\mathbb{R}^{d}$ (ie. $f$ is
identified to a set-valued map $F$), we have:\smallskip\newline(i).
$DF(\bar{x},f(\bar{x}))(u)\neq\emptyset,$ for every $u\in\mathbb{R}^{d}$ ; and
\smallskip\newline(ii) for all $(u,v)\in\mathrm{gph}(DF(\bar{x},f(\bar{x})))$
it holds: $|v|\,\leq k_{\bar{x}}\,|u|.$
\end{lemma}

\noindent\textbf{Proof.} Let $u\in\mathbb{R}^{d}$ and $(h_{n})\searrow0^{+}.$
We deduce from \eqref{eq:calm} that for $n$ sufficiently large:
\[
w_{n}:=\frac{f(\bar{x}+h_{n}u)-f(\bar{x})}{h_{n}}\in(1+k_{\bar{x}})\mathrm{B}.
\]
Since the sequence $(w_{n})_{n}$ is bounded, passing to a subsequence we may
assume that it converges to some $v\in\mathbb{R}^{\ell}$. Since $f(\bar
{x})+h_{n}w_{n}=f(\bar{x}+h_{n}u)\in F(\bar{x}+h_{n}u)$ we conclude
from~\eqref{DFk} that $v\in DF(\bar{x},f(\bar{x}))(u)$ and (i)
holds.\smallskip\newline To prove (ii), fix $(u,v)\in\mathrm{gph}(DF(\bar
{x},f(\bar{x})))$ and consider the associated sequences $(h_{n})\searrow
0^{+}\,$and $(u_{n},v_{n})_{n}\rightarrow(u,v)$ such that $f(\bar{x})+h_{n}v_{n}=f(\bar{x}+h_{n}u_{n})$. In view of~\eqref{eq:calm} there exists a
sequence of positive numbers $\varepsilon_{n}\rightarrow0^{+}$ such that
\[
\displaystyle|v_{n}|=\frac{|f(\bar{x}+h_{n}u_{n})-f(\bar{x})|}{h_{n}}
\leq(k_{\bar{x}}+\varepsilon_{n})|u_{n}|,\;\forall n\geq0.
\]
Passing to the limit we obtain $|v|\leq k_{\bar{x}}\,|u|$ as desired.\hfill
$\Box$

\bigskip

\noindent Notice that every Lipschitz function at $\bar{x}$ satisfies~\eqref{eq:calm}. Notice that in Lemma~\ref{Lemma-a}, $f$ is
considered both as a (single-valued)\ function and as a set-valued map with singleton values, and consequently we can define the usual (Fr\'echet) derivative and the graphical derivative. The following proposition clarifies the relation between these two objects.

\begin{proposition}
[compatibility]\label{df}Let $f:\mathbb{R}^{d}\longrightarrow\mathbb{R}^{\ell
}$ be a continuous function satisfying~\eqref{eq:calm} at $\bar{x}\in\mathbb{R}^{n}.$ Then $f$ is differentiable at $\bar{x}$ if and only if the map $F(x)=\{f(x)\}$ is differentiable at $(\bar{x},f(\bar{x}))$. In this case,
it holds
\[
\mathrm{gph}(DF(\bar{x},f(\bar{x})))=\mathrm{gph}(df(\bar{x})).
\]

\end{proposition}

\noindent\textbf{Proof.} Let us first assume that $f$ is differentiable at
$\bar{x},$ that is,
\begin{equation}
\displaystyle\lim_{|z|\rightarrow0}\frac{|f(\bar{x}+z)-f(\bar{x})-df(\bar
{x})(z)|}{|z|}=0. \label{diff}
\end{equation}
Then it follows easily from~\eqref{DFk} that for any $u\in\mathbb{R}^{d}$ it
holds $df(\bar{x})(u)\in DF(\bar{x},\bar{y})(u)$, therefore $\mathrm{gph}(df(\bar{x}))\subset\mathrm{gph}(DF(\bar{x},f(\bar{x}))).$ Let now $v\in
DF(\bar{x},f(\bar{x}))(u)$. It follows from~\eqref{DF} that there exist
sequences $(h_{n})\searrow0^{+}$ and $(u_{n},v_{n})_{n}\rightarrow(u,v)$ such
that$\;f(\bar{x})+h_{n}v_{n}\in F(\bar{x}+h_{n}u_{n})=\{f(\bar{x}+h_{n}u_{n})\},$ for all $n\geq1$. It follows readily that
\[
\displaystyle v_{n}=\frac{f(\bar{x}+h_{n}u_{n})-f(\bar{x})}{h_{n}}.
\]
In view of \eqref{diff}, taking the limit as $n\rightarrow\infty$ we obtain
$v=df(\bar{x})(u)$, yielding the equality
\[
\mathrm{gph}(DF(\bar{x},f(\bar{x})))=\mathrm{gph}\left(  df(\bar{x})\right)
.
\]
The above ensures the differentiability of $F$ at $(\bar{x},f(\bar{x}))$,
since $\mathrm{gph}\left(  df(\bar{x})\right)  $ is a linear subspace of
$\mathbb{R}^{d}\times\mathbb{R}^{\ell}$ and consequently a closed convex
cone.\medskip

Let us now assume that $F$ is differentiable at $(\bar{x},f(\bar{x}))$, that
is, $DF(\bar{x},f(\bar{x}))$ is a convex process and consequently
$\mathrm{gph}\left(  DF(\bar{x},f(\bar{x}))\right)  $ is a closed convex cone.
We are going to prove successively that $\mathrm{gph}\left(  DF(\bar{x},f(\bar{x}))\right)  $ is the graph of some linear function $T$ from
$\mathbb{R}^{d}$ to $\mathbb{R}^{\ell}$ and then that $T$ is the Fr\'echet
differential of $f$ at $\bar{x}$.\smallskip\newline To this end, fix
$(u,v)\in\mathrm{gph}(DF(\bar{x},f(\bar{x})))$. We claim that
\begin{equation}
DF(\bar{x},f(\bar{x}))(-u)=\{-v\}. \label{eq:m1}
\end{equation}
In view of Lemma~\ref{Lemma-a}(ii), this is clear when $u=0$, therefore, we
may assume $u\neq0$. Let $w\in DF(\bar{x},f(\bar{x}))(-u)$ and notice that
$(u,v)$ and $(-u,w)$ are two points of the convex set $\mathrm{gph}\left(
DF(\bar{x},f(\bar{x}))\right)  $. It follows that
\[
\frac{1}{2}[(u,v)+(-u,w)]=(0,\frac{1}{2}(v+w)\in DF(\bar{x},f(\bar{x})),
\]
which again in view of Lemma~\ref{Lemma-a}(ii) yields $v+w=0$
and~\eqref{eq:m1} follows. Moreover, since this is true for any $u\in
\mathbb{R}^{d}$, we conclude that the set $DF(\bar{x},f(\bar{x}))(u)$ reduces
to the singleton $\{v\}$. We claim that for every $u$ the directional
derivative
\begin{equation}
\displaystyle\lim_{t\rightarrow0}\,\frac{f(\bar{x}+tu)-f(\bar{x})}{t}
\label{eq:m2}
\end{equation}
exists and is equal to $v$. Indeed let $t_{n}\searrow0^{+}$ and observe that
for $n$ sufficiently large the sequence
\[
v_{n}:=\frac{f(\bar{x}+t_{n}u)-f(\bar{x})}{t_{n}}
\]
is contained to the ball $B(0,$ $1+k_{\bar{x}})$ therefore it has accumulation
points, which all belong to $DF(\bar{x},f(\bar{x}))(u)=\{v\}$. This proves
that $v_{n}\rightarrow v$ and that the directional limit in~\eqref{eq:m2}
exists when $t\rightarrow0^{+}$ and is equal to $DF(\bar{x},f(\bar{x}))(u)$.
Changing $u$ into $-u$ we obtain that the corresponding directional limit also
exists and it is equal to $-Df(\bar{x},f(\bar{x}))(-u),$ which in turn, in
view of~\eqref{eq:m1}, is equal to $DF(\bar{x},f(\bar{x}))(u)$. It remains to
show that the linear function
\[
\mathbb{R}^{d}\ni u\mapsto T(u):=v\in\mathbb{R}^{d}\qquad\text{where }
DF(\bar{x},f(\bar{x}))(u)=\{v\},
\]
is the derivative of $f$ at $\bar{x}$ namely:
\[
\displaystyle\lim_{|z|\rightarrow0^{+}}\frac{f(\bar{x}+z)-f(\bar{x})-T(z)}{|z|}=0.
\]
Let us assume, towards a contradiction, that this not the case. Then there
exist $\alpha>0$ and a sequence $(z_{n})\rightarrow0$, with $z_{n}\neq0$ for
all $n\geq1$ such that
\begin{equation}
\frac{|f(\bar{x}+z_{n})-f(\bar{x})-T(z_{n})|}{|z_{n}|}\geq\alpha,\;\forall
n\geq0. \label{www}
\end{equation}
The sequence $u_{n}:=z_{n}/|z_{n}|$ is bounded and converges (up to a
subsequence) to some $u\in\mathbb{R}^{d}$. Moreover, using~\eqref{eq:calm} for
$n$ sufficiently large we obtain
\[
\begin{split}
\left\vert \frac{f(\bar{x}+z_{n})-f(\bar{x})-T(z_{n})}{z_{n}}\right\vert  &
=\left\vert \frac{f(\bar{x}+|z_{n}|u_{n})-f(\bar{x})}{|z_{n}|}-T(u_{n})\right\vert \leq\left\vert \frac{f(\bar{x}+|z_{n}|u)-f(\bar{x})}{|z_{n}%
|}-T(u)\right\vert \\
&  +\,\frac{|f(\bar{x}+|z_{n}|u_{n})-f(\bar{x}+|z_{n}|u)|}{|z_{n}|}+|T(u)-T(u_{n})|\\
&  \leq\,\left\vert \frac{f(\bar{x}+|z_{n}|u)-f(\bar{x})}{|z_{n}|}-T(u)\right\vert \,+\,(1+k_{\bar{x}})\,|u_{n}-u|\,+\Vert T\Vert\,|u-u_{n}|
\end{split}
\]
where $\Vert T\Vert$ denotes the operator norm of the linear map $T$. Passing
to the limsup to the above inequalities and taking into account that
\[
u_{n}\rightarrow u\qquad\text{and \qquad}\lim_{n}\,\frac{f(\bar{x}+|z_{n}|u)-f(\bar{x})}{|z_{n}|}=T(u)
\]
we obtain
\[
\limsup_{n}\frac{|f(\bar{x}+z_{n})-f(\bar{x})-T(z_{n})|}{|z_{n}|}\leq0
\]
which contradicts \eqref{www}. Therefore, $T(u)=df(\bar{x})(u)$ and
$DF(\bar{x},f(\bar{x}))(u)=\{df(\bar{x})(u)\}$. \\
The proof is complete.
\hfill$\Box$

\bigskip

\section{Set-valued maps generated by a finite family}

\label{sec:4}In this section we consider a particular type of set-valued maps
$F:\mathbb{R}^{d}\rightrightarrows\mathbb{R}^{\ell}$, namely those generated
by a finite family of locally Lipschitz functions, that is,
\[
F(x)=\mathrm{conv}\{f_{1}(x),\ldots,f_{N}(x)\}
\]
where $f_{i}:\mathbb{R}^{d}\rightarrow\mathbb{R}^{\ell}$ is a (locally)
Lipschitz function, for all $i\in\{1,\ldots,N\}$. It is easy to see that
$F:\mathbb{R}^{d}\rightrightarrows\mathbb{R}^{\ell}$ is Lipschitz continuous
(as set-valued function) with nonempty convex compact values. We shall further
assume that for every $x\in\mathbb{R}^{d}$ the family $\{f_{1}(x),\ldots
,f_{N}(x)\}$ is affinely independent. This assumption yields in particular
that $N\leq\ell+1.$ Moreover, setting
\[
\Delta:=\left\{  (\lambda_{1},\ldots,\lambda_{N})\in\lbrack0,1]^{N}:\;
{\displaystyle\sum\limits_{j=1}^{N}}
\lambda_{j}=1\right\}  ,
\]
we have a unique representation of every $y\in F(x),$ that is, there exist
unique $(\lambda_{1},\ldots,\lambda_{N})\in\Delta$ such that $y=
{\displaystyle\sum\limits_{j=1}^{N}}
\lambda_{j}f_{j}(x)$.\smallskip\newline
We are now ready to state the following generalization of Rademacher theorem for this type of set-valued maps.
\begin{proposition} [generalized Rademacher theorem]\label{cor:family} Let $f_{i}:\mathbb{R}^{d}\rightarrow\mathbb{R}^{\ell}$, $i\in\{1,\dots, N\}$ be a family of Lipschitz continuous functions for which $\{f_{1}(x),\ldots,f_{N}(x)\}$ are affinely independent for
every $x\in\mathbb{R}^{d}$.  Then the set-valued map $F:\mathbb{R}^{d}\rightrightarrows\mathbb{R}^{\ell}$ given by
\[
F(x)=\mathrm{conv\,}\{f_{1}(x),\ldots,f_{N}(x)\},\quad x\in\mathbb{R}^{d}\text{,}
\]
is differentiable at $(\bar x,\bar y)$ for \emph{a.e.} $\bar x\in\mathbb{R}^{d}$ and all $\bar y\in F(x)$. In particular, for every $u\in\mathbb{R}^{d}$
\[
DF(\bar{x},\bar{y})(u)\,=\,
{\displaystyle\sum\limits_{j=1}^{N}}
\bar{\lambda}_{j}\,df_{j}(\bar{x})(u)\,+\,T_{F(\bar{x})}(\bar{y}),
\]
where $(\bar{\lambda}_{1},\ldots,\bar{\lambda}_{N})\in\Delta$ are such that
\begin{equation}
\bar{y}={\displaystyle\sum\limits_{j=1}^{N}}
\bar{\lambda}_{j}\,f_{j}(\bar{x}). \label{eq:a1}
\end{equation}

\end{proposition}

\noindent\textbf{Proof.} Let us denote by $\mathcal{D}\subset\mathbb{R}^{d}$ the set of points of common
differentiability of the functions $f_{i}$, $i\in\{1,\dots,N\}$ (by Rademacher theorem this set is of full measure) and let us fix $\bar{x}\in\mathcal{D}$,
$\bar{y}\in F(\bar{x})$ and $u\in\mathbb{R}^{d}$ . Let further $w\in DF(\bar{x},\bar{y})(u)$. Then
there exist $(w_{n})_{n}\rightarrow w$ and $(t_{n})\searrow0^{+}$ such that
\[
\bar{y}+t_{n}w_{n}\in F(\bar{x}+t_{n}u).
\]
Then for every $n\geq1$ there exists $(\lambda_{1}^{n},\ldots,\lambda_{N}^{n})\in\Delta$ such that
\begin{equation}
\bar{y}+t_{n}w_{n}=
{\displaystyle\sum\limits_{j=1}^{N}}
\lambda_{j}^{n}\,f_{j}(\bar{x}+t_{n}u). \label{eq:a0}
\end{equation}
Thanks to the compactness of $\Delta$, passing to a subsequence, we
may assume that $$(\lambda_{1}^{n},\ldots,\lambda_{N}^{n})\rightarrow
(\lambda_{1}^{\ast},\ldots,\lambda_{N}^{\ast}),\quad \text{as} n\rightarrow\infty.$$ For
every $j\in\{1,\ldots,N\}$, the first-order Taylor series of $f_{j}$ at
$\bar{x}$ yields:
\begin{equation}
f_{j}(\bar{x}+t_{n}u)=f_{j}(\bar{x})\,+\,t_{n}\,df_{j}(\bar{x})(u)\,+o_{j}
(t_{n}u),\qquad\text{where\ }\lim_{|v|\rightarrow0}\,\frac{o_{j}(v)}{|v|}=0.
\label{eq:ar11}
\end{equation}
Replacing~\eqref{eq:a1}, \eqref{eq:ar11} into~\eqref{eq:a0} we obtain:
\begin{equation}
t_{n}w_{n}\,=\,
{\displaystyle\sum\limits_{j=1}^{N}}
\left(  \lambda_{j}^{n}-\bar{\lambda}_{j}\right)  f_{j}(\bar{x})\,+\,t_{n}
{\displaystyle\sum\limits_{j=1}^{N}}
\lambda_{j}^{n}\,df_{j}(\bar{x})(u)\,+\,t_{n}\,||u||\,
{\displaystyle\sum\limits_{_{n}j=1}^{N}}
\lambda_{j}^{n}\left(  \frac{o_{j}(t_{n}u)}{||t_{n}u||}\right)  ,
\label{eq:a2}
\end{equation}
and passing to the limit as $n\rightarrow\infty$ we deduce:
\[
{\displaystyle\sum\limits_{j=1}^{N}}
\left(  \lambda_{j}^{\ast}-\bar{\lambda}_{j}\right)  f_{j}(\bar{x}
)=0,\qquad\text{yielding}\quad\lambda_{j}^{\ast}=\bar{\lambda}_{j},\text{ for
all }j\in\{1,\ldots,N\}.
\]
Returning to~(\ref{eq:a0}), we deduce easily that:
\begin{equation}
w_{n}\,=\,
{\displaystyle\sum\limits_{j=1}^{N}}
\lambda_{j}^{n}\left(  \frac{f_{j}(\bar{x}+t_{n}u)-f_{j}(\bar{x})}{t_{n}
}\right)  \,+\,\frac{1}{t_{n}}\,\left(
{\displaystyle\sum\limits_{j=1}^{N}}
\,\lambda_{j}^{n}\,f_{j}(\bar{x})-\bar{y}\right)  . \label{eq:a3}
\end{equation}
Noticing that
\[
\frac{1}{t_{n}}\left(
{\displaystyle\sum\limits_{j=1}^{N}}
\,\lambda_{j}^{n}\,f_{j}(\bar{x})-\bar{y}\right)  \,\in\,\mathbb{R}_{+}\left(
F(\bar{x})-\bar{y}\right)  \subset T_{F(\bar{x})}(\bar{y})
\]
and taking the limit in~(\ref{eq:a3}) as $n\rightarrow\infty$ we deduce that
\[
w=\,\lim_{n\rightarrow\infty}\,w_{n}\,=\,
{\displaystyle\sum\limits_{j=1}^{N}}
\bar{\lambda}_{j}df_{j}(\bar{x})(u)+T_{F(\bar{x})}(\bar{y}).
\]

For the reverse inclusion, it suffices to prove that for every $u\in
\mathbb{R}^{d}$ we have:
\[
{\displaystyle\sum\limits_{j=1}^{N}}
\bar{\lambda}_{j}df_{j}(\bar{x})(u)+\mathbb{R}_{+}\left(  F(\bar{x})-\bar
{y}\right)  \subset DF(\bar{x},\bar{y})(u).
\]
To this end, fix $\mu\geq0,$ $(\lambda_{1},\ldots,\lambda_{N})\in\Delta$ and
define
\[
v=\mu\left(
{\displaystyle\sum\limits_{k=1}^{N}}
\,\lambda_{k}\,f_{k}(\bar{x})-\bar{y}\right)  \qquad\text{(arbitrary element
of }\mathbb{R}_{+}\left(  F(\bar{x})-\bar{y}\right)  \text{).}
\]
Let further $(t_{n})\searrow0^{+}$ and set
\begin{equation}
w_{n}:=\frac{1}{t_{n}}\left(  (1-\mu t_{n})
{\displaystyle\sum\limits_{j=1}^{N}}
\,\bar{\lambda}_{j}\,f_{j}(\bar{x}+t_{n}u)+\mu t_{n}
{\displaystyle\sum\limits_{k=1}^{N}}
\,\lambda_{k}\,f_{k}(\bar{x}+t_{n}u)-\bar{y}\right)  . \label{eq:ar1}
\end{equation}
Since $F(\bar{x}+t_{n}u)$ is convex, the above definition yields
\[
\bar{y}+t_{n}w_{n}\in F(\bar{x}+t_{n}u),\quad\text{for all }n\text{sufficiently large,}
\]
and consequently
\[
w=\,\lim_{n\rightarrow\infty}\,w_{n}\,\in DF(\bar{x},\bar{y})(u).
\]
Replacing $\bar{y}$ in~(\ref{eq:ar1}) by its representation given in
(\ref{eq:a1}) we deduce:
\[
w_{n}:=(1-\mu t_{n})
{\displaystyle\sum\limits_{j=1}^{N}}
\,\bar{\lambda}_{j}\,\frac{f_{j}(\bar{x}+t_{n}u)-f_{j}(\bar{x})}{t_{n}}
+\mu\left(
{\displaystyle\sum\limits_{k=1}^{N}}
\,\lambda_{k}\,f_{k}(\bar{x}+t_{n}u)-\bar{y}\right)  .
\]
Taking the limit as $n\rightarrow\infty$ we obtain:
\[
{\displaystyle\sum\limits_{j=1}^{N}}
\,\bar{\lambda}_{j}\,df_{j}(\bar{x})(u)+v\,\in DF(\bar{x},\bar{y})(u),
\]
as desired. The proof is complete. \hfill$\Box$

\bigskip

\begin{remark}
Taking $N=1$ in Proposition~\ref{cor:family} we recover the classical Rademacher theorem.
\end{remark}


\section{A Rademacher Theorem for isotropically Lipschitz maps}\label{sec:5}

In this chapter we deal with general set-valued maps for which we do not dispose a concrete description. Let us first notice that if a Lipschitz set-valued map $F:\mathbb{R}^d\rightrightarrows \mathbb{R}^{\ell}$ has closed convex values with nonempty interior, then differentiability of $F$ at $(x,y)$ for every $x\in\mathbb{R}^d$ and for \emph{a.e.} $y\in F(x)$ (with respect to the $\mathcal{L}_{\ell}$-measure of $F(x)$) is straightforward.  Indeed, for every $y\in\mathrm{int} F(x)$ we have $\mathrm{ghp} \left(DF(x,y)\right)=\mathbb{R}^d\times\mathbb{R}^{\ell}$ and the boundary $\partial F(x)$ is a null subset of the convex set $F(x)$. In this case, the interesting part is clearly the behavior of $F$ around boundary points $y\in\partial F(x)$ which is what we ar going to investigate in the forthcoming Theorem~\ref{RadeFint}. \smallskip\newline
Let $F:\mathbb{R}^{d}\rightrightarrows\mathbb{R}^{\ell}$ be a
Lipschitz set-valued map with nonempty compact values. For every
$p\in\mathbb{R}^{\ell}$ we set:
\[
\sigma_{F(x)}(p):=\max_{y\in F(x)}\,\langle p,y\rangle\quad\text{and\quad
}Y(x,p):=\{y\in F(x):\,\langle p,y\rangle=\sigma_{F(x)}(p)\}.
\]
Notice that $Y(x,p)$ is a nonempty compact subset of $\mathbb{R}^{\ell}$ and that the function $p\mapsto\sigma_{F(x)}(p)$ is convex and
Lipschitz continuous. It is also easy to verify that
\[
(x,p)\mapsto\sigma_{F(x)}(p)\quad\text{is (locally) Lipschitz continuous,}
\]
while the set-valued function $(x,p)\rightrightarrows Y(x,p)$ is continuous,
but in general not Lipschitz continuous, as shows the following example.

\begin{example} \label{contrexemple} Consider the real valued functions
\[
\label{obs2}
\begin{array}
[c]{lll}
\phi_{0}(x) & := & \left\{
\begin{array}
[c]{lll}
x^{2}-x & \mbox{ if } & x\leq0\\
x^{2} & \mbox{ if } & x>0,
\end{array}
\right.
\end{array}
\]
and for $\tau>0$:
\[
\label{obs3}
\begin{array}
[c]{lll}
\phi_{\tau}(x) & := & \left\{
\begin{array}
[c]{lll}
x^{2}-x+\tau & \mbox{ if } & x\leq0\\
-x\sqrt{\tau}+\tau & \mbox{ if } & x\in\lbrack0,\sqrt{\tau}],\\
x^{2}-\tau & \mbox{ if } & x>\sqrt{\tau},
\end{array}
\right.
\end{array}
\]
We define $F:[-1,1]\rightrightarrows\mathbb{R}^{2}$ by
\[
F(\tau)=\mathrm{epi\,}\phi_{|\tau|}\cap(\mathbb{R}\times\lbrack-2,2]).
\]
Then $F$ is Lipschitz continuous with compact convex values (with a slight
modification we can obtain that $F(\tau)$ is strictly convex, but for $\bar p =(0,1)$,
the (single-valued) map $\tau\rightarrow Y(\tau,\bar p):=\{(\sqrt{\tau}, 0)\}$  is not Lipschitz continuous.
\end{example}

\noindent The above example motivates the following definition.

\begin{definition}
[isotropically Lipschitz map]\label{def_isotropic}Let $F:\mathbb{R}^{d}\rightrightarrows\mathbb{R}^{\ell}$ be a set-valued map with nonempty
convex compact values. We say that $F$ is (locally) isotropically Lipschitz
continuous, if (for every $\bar{x}\in\mathbb{R}^{d}$) there exists $k>0$ such
that for all $p\in\mathbb{R}^{\ell}$ and all $x,x^{\prime}$ (in a neighborhood
of $\bar{x}$) it holds:
\begin{equation}
Y(x,p)\subset Y(x^{\prime},p)+k\,|x-x^{\prime}|\,\mathrm{B.}
\label{eq:isotropic}
\end{equation}

\end{definition}

\noindent The following result holds.

\begin{proposition}
[isotropically Lipschitz vs Lipschitz]\label{Prop_isotropic}Every (locally)
isotropically Lipschitz map $F:\mathbb{R}^{d}\rightrightarrows\mathbb{R}^{\ell}$ is (locally) Lipschitz.
\end{proposition}

\noindent\textbf{Proof.} Let $F$ be $k$-isotropically Lipschitz around $\bar{x}\in\mathbb{R}^{d}$
and let us assume, towards a contradiction, that $F$ is not $k$-Lipschitz there. Then there exist sequences
$(x_{n})_{n}$, $(x_{n}^{\prime})_{n}$ converging to $\bar{x}$ and $z_{n}\in F(x_{n})$ such that
$z_{n}\notin F(x_{n}^{\prime})+\,k\,|x_{n}-x_{n}^{\prime
}|\,\mathrm{B.}$ Denoting by $z_{n}^{\prime}$ the projection of $z_{n}$ onto
the convex compact set $F(x_{n}^{\prime})$ and setting
\[
p_{n}=\frac{z_{n}-z_{n}^{\prime}}{|z_{n}-z_{n}^{\prime}|}
\]
we deduce that $z_{n}^{\prime}\in Y(x_{n}^{\prime},p_{n})$ and
\[
\langle p_{n},z_{n}-z_{n}^{\prime}\rangle=|z_{n}-z_{n}^{\prime}|\geq
k\,|x_{n}-x_{n}^{\prime}|.
\]
Then for every $y_{n}\in Y(x_{n},p_{n})\subset F(x_{n})$ one has:
\[
|y_{n}-z_{n}^{\prime}|\,\geq\,\langle p_{n},y_{n}-z_{n}^{\prime}\rangle
\,\geq\,\langle p_{n},z_{n}-z_{n}^{\prime}\rangle=|z_{n}-z_{n}^{\prime}|\,\geq
k\,|x_{n}-x_{n}^{\prime}|.
\]
The above clearly yields
\[
z_{n}^{\prime}\notin Y(x_{n},p_{n})+\,k\,||x_{n}-x_{n}^{\prime
}|\,\mathrm{B}
\]
contradicting (\ref{eq:isotropic}).\hfill$\Box$

\bigskip

\begin{remark}
$\mathrm{(i).}$ If $F$ is a (single-valued) function, then the notions of isotropic
Lipschitz continuity and Lipschitz continuity coincide.\smallskip
\newline$\mathrm{(ii).}$ If $F\ $is isotropically Lipschitz and has strictly convex
values, then for every $p\in S^{\ell-1}$ the set-valued map
\[
x\rightrightarrows\lbrack Y(x,-p),Y(x,p)]
\]
is Lipschitz. (Notice that strict convexity of $F(x)$ guarantees that both
$Y(x,p)$ and $Y(x,-p)$ are singletons, so that $[Y(x,-p),Y(x,p)]$ is a closed
segment in $\mathbb{R}^{\ell}$.)
\end{remark}

\medskip

\noindent We shall show that locally isotropically Lipschitz set-valued maps
having \textit{strictly convex} values with \textit{nonempty interior} satisfy
a Rademacher-type result. Before going further, let us make some comments
about this last assumption.\smallskip\newline(A). Fix $\bar{x}\in
\mathbb{R}^{d}$ and set $K=F(\bar{x})\subset\mathbb{R}^{\ell}.$ Assuming that
the strictly convex compact set $K$ has nonempty interior guarantees that its
boundary $\partial K$ is locally the graph of a (strictly) convex Lipschitz function
$g:\mathbb{R}^{\ell-1}\mapsto\mathbb{R}$, that is, for every $y_{0}\in K$
there exists a (strictly) convex Lipschitz function $g$ and $\varepsilon>0$ such that
\[
\partial K\cap B(y_{0},\varepsilon)=\mathrm{gph}(g)\cap B(y_{0},\varepsilon).
\]
It follows that every $y\in\partial K\cap B(y_{0},\varepsilon)$ can be
represented in local coordinates as $y=(\xi,g(\xi))$ and that for
\textit{a.e.} $\xi\in\mathbb{R}^{\ell-1}$ (with respect to the Lebesgue
measure of $\mathbb{R}^{\ell-1}$) the gradient $\nabla g(\xi)$ exists and
therefore%
\[
N_{K}((\xi,g(\xi))=\mathbb{R}_{+}(\nabla g(\xi),-1).
\]
Since both the projection $\pi(y)=\pi(\xi,g(\xi))=\xi\in\mathbb{R}^{\ell-1}$
and its inverse $\pi^{-1}(\xi)=(\xi,g(\xi))$ are Lipschitz continuous, we get
a bi-Lipschitz bijection between $\partial K\cap B(y_{0},\varepsilon)$ and
$\mathrm{gph}(g)\cap B(y_{0},\varepsilon)$ and conclude that these sets have
the same null subsets (in the sense that they are identified through this bijection). Therefore for \emph{a.e} $\bar{y}\in\partial F(\bar{x})$
we have $N_{\partial F(\bar{x})}(\bar{y})=\mathbb{R}\bar{p}$ and $\bar
{y}=Y(\bar{x},\bar{p})$, which means that $\bar{y}$ is a point of smoothness
of the boundary of $\partial F(\bar{x})$ and
\begin{equation}
\mathcal{M}_{\bar{x}}:=\{\bar{y}\in\partial F(\bar{x}):\;\mbox{
dim }N_{\partial F(\bar{x})}(\bar{y})=1\,\}\label{eq:Mx}
\end{equation}
has a full measure in $\partial F(\bar{x})$.\medskip\newline(B). Thanks to the
assumption of  strict convexity of $F(\bar{x}),$ the mapping
\[
y_{\bar{x}}:p\in S^{\ell-1}\mapsto Y(\bar{x},p)
\]
is single-valued. It is also clearly surjective (by Hahn-Banach theorem). Set
\[
\mathbf{M}_{\bar{x}}:=\{p\in S^{\ell-1}:\;y_{\bar{x}}^{-1}(y_{\bar{x}}(p))=\{p\}\}=\{p\in S^{\ell-1}:\;y_{\bar{x}}(p)\in\mathcal{M}_{\bar{x}}\}.
\]
Then $p\mapsto y_{\bar{x}}$ is a bijection from $\mathbf{M}_{\bar{x}}$ onto
$\mathcal{M}_{\bar{x}}$ and $N_{\partial F(\bar{x})}(y_{\bar{x}}(p))=\mathbb{R}_{+}p$ for all $p\in\mathcal{M}_{\bar{x}}$. Notice further that
$\mathcal{M}_{\bar{x}}$ is full-measure in $\partial F(\bar{x})$ and
$\mathbf{M}_{\bar{x}}$ is full-measure in $S^{\ell-1}$. Moreover, for any
full-measure subset $P$ of $S^{\ell-1}$, the set $y_{\bar{x}}(S^{\ell
-1}\backslash P)$ is null in $\partial F(\bar{x})$ and $\mathcal{M}_{\bar{x}}\diagdown y_{\bar{x}}(S^{\ell-1}\backslash P)$ is full-measure in $\partial
F(\bar{x})$.

\bigskip

\noindent We are now ready to state a Rademacher type result for
isotropically Lipschitz set-valued maps.

\begin{theorem}
[Rademacher result for isotropically Lipschitz maps.]\label{RadeFint} Let
$F:\mathbb{R}^{d}\rightrightarrows\mathbb{R}^{\ell}$ be a locally
isotropically Lipschitz set-valued map, with strictly convex values of
nonempty interior. Then for \emph{a.e.} $\bar{x}\in\mathbb{R}^{d}$ $($with respect to the $\mathcal{L}_d$-Lebesgue measure of $\mathbb{R}^d)$ and for \emph{a.e.} $\bar{y}\in\partial F(\bar{x})$ $($with respect to the $(\ell-1)$-dimensional measure of $\partial F(x))$, $F$ is differentiable at $(\bar x, \bar y)$
and moreover,
\[
\mathrm{gph}(DF(\bar{x},\bar{y}))\text{\quad is a half-space.}
\]

\end{theorem}

\noindent\textbf{Proof.} The mapping $(x,p)\mapsto\sigma_{F(x)}(p)$ is
Lipschitz. By Rademacher theorem, there exists a set $\hat{D}\subset
\mathbb{R}^{d}\times\mathbb{R}^{\ell}$ of full measure such that for any
$(x,p)\in\hat{D}$ the derivative $\nabla\sigma_{F(x)}(p)$ exists. By Fubini
theorem, and using the homogeneity of $\sigma$ with respect to $p$ we deduce
the existence of a full measure subset $D\subset\mathbb{R}^{d}$ such that for
any $x\in D$ the set
\[
\{p\in S^{\ell-1}:\;(x,p)\in\hat{D}\}\;\text{has a full measure in }S^{\ell
-1}.
\]
Fix $\bar{x}\in D$ and set
\[
P_{\bar{x}}:=\{p\in S^{\ell-1}:\;(x,p)\in\hat{D}\}\quad\text{(which has
full-measure in }S^{\ell-1}\text{).}
\]
Recalling the definition of $\mathcal{M}_{\bar{x}}$ from (\ref{eq:Mx}) we set
\[
\widehat{\mathcal{M}_{\bar{x}}}=\mathcal{M}_{\bar{x}}\diagdown y_{\bar{x}
}(S^{\ell-1}\backslash P_{\bar{x}})\quad\text{(which has full-measure in
}\partial F(\bar{x})\text{).}
\]
Fix $\bar{y}\in\widehat{\mathcal{M}_{\bar{x}}}$ and set $\bar{p}=y_{\bar{x}
}^{-1}(\bar{y})$ (i.e. $Y(\bar{x},\bar{p})=\{\bar{y}\}$). Set
\[
C_{\bar{x},\bar{y}}:=\{\,(u,w)\in\mathbb{R}^{d}\times\mathbb{R}^{\ell
}:\;\nabla_{x}\sigma_{F(\bar{x})}(\bar{p})(u)\geq\langle\bar{p},w\rangle\,\}
\]
which is a closed half-space. We are going to prove that
\begin{equation}
T_{\mathrm{gph}(F)}(\bar{x},\bar{y})=C_{\bar{x},\bar{y}}\label{TP}
\end{equation}
which immediately yields that the graphical derivative is convex (a half-space
in $\mathbb{R}^{d}\times\mathbb{R}^{\ell}$) and consequently $F$ is
differentiable at $(\bar{x},\bar{y})$, according to Definition~\ref{diffF}.\smallskip\newline To this end, take $(u,w)\in T_{\mathrm{gph}(F)}(\bar
{x},\bar{y})$. In view of Proposition~\ref{Prop_isotropic} and (\ref{DFk}),
there exist sequences $t_{n}\searrow0^{+}$ and $w_{n}\rightarrow w$ such that
$\bar{y}\in F(\bar{x}+t_{n}u)-t_{n}w_{n}.$ \smallskip\newline
For $\bar{p}=y_{\bar{x}}^{-1}(\bar{y})$, we have $Y(\bar{x},\bar
{p})=\{\bar{y}\}$ and
\[
\sigma_{F(\bar{x})}(\bar{p})=\langle\bar{p},\bar{y}\rangle\,\leq\,\sigma
_{F(\bar{x}+t_{n}u)}(\bar{p})-t_{n}\langle\bar{p},w_{n}\rangle,\quad\text{yielding}\quad
\frac{\sigma_{F(\bar{x})+t_{n}w_{n}}(\bar{p})-\sigma_{F(\bar{x}
)}(\bar{p})}{t_n}\,\geq\,\langle w_{n},\bar{p}\rangle.
\]
Therefore, $\nabla_{x}\sigma_{F(\bar{x})}(\bar
{p})(u)\geq\langle\bar{p},w\rangle$ and the inclusion
$T_{\mathrm{gph}(F)}(\bar{x},\bar{y})\subset C_{\bar{x},\bar{y}}$ follows.\smallskip\newline
For the opposite inclusion, since
$T_{\mathrm{gph}(F)}(\bar{x},\bar{y})$ is closed cone, it suffices to show
that
\[
\mathrm{int\,}C_{\bar{x},\bar{y}}\subset T_{\mathrm{gph}(F)}(\bar{x},\bar{y}).
\]
To this end, let $(u,w)\in\mathbb{R}^{d}\times\mathbb{R}^{\ell}$ be such that
\begin{equation}\label{eq:adan}
\nabla_{x}\sigma_{F(\bar{x})}(\bar{p})(u)\,>\,\langle\bar{p},w\rangle.
\end{equation}
We may also assume that $w\notin\mathbb{R}\bar{p}$ (if $w\in\mathbb{R}\bar{p}$
then we replace $w$ by another vector $\tilde{w}\notin\mathbb{R}\bar{p}$
arbitrarily close to it). Let us assume by contradiction that $(u,w)\notin
T_{\mathrm{gph}(F)}(\bar{x},\bar{y})$. Then by a standard argument (see \cite{Q1992} \emph{e.g.})
\begin{equation}
\exists\delta>0,\;\forall\tau\in(0,\delta):\quad(\bar{y}+\tau w+\tau
\delta\mathrm{B})\cap F(\bar{x}+\tau u)=\emptyset.\label{contr}
\end{equation}
Let us first prove that $u\neq0$. Indeed if $u=0,$ then $\nabla_{x}\sigma_{F(\bar
{x})}(\bar{p})(u)=0>\,\langle\bar{p},w\rangle$ and $(\bar{y}+\tau w+\tau
\delta\mathrm{B})\cap F(\bar{x})=\emptyset$ for all $\tau\in(0,\delta)$.
Therefore $w\notin T_{F(\bar{x})}(\bar{y}),$ which is the half-space $[\bar
{p}\leq0]$ because $\partial F(\bar{x})$ is smooth at $\bar{y}$ and
$\mathrm{int\,}F(\bar{x})\neq\emptyset$. This is a contradiction with
$\langle\bar{p},w\rangle<0$.\smallskip\newline
Set
\begin{equation}\label{eq:ada}
\left\{
\begin{array}
[c]{l}
e_{1}:=w-\langle\bar{p},w\rangle\,\bar{p}\\
e_{2}:=-\bar{p}
\end{array}
\right.  \text{\qquad and\qquad}Z:=\mathrm{span}\{e_{1},e_{2}\}.
\end{equation}
We deduce from \eqref{contr} that
\begin{equation}
\Bigl [ \{(\bar{y}+\tau w+\tau\delta\mathrm{B})-Y(\bar{x}+\tau u,\bar
{p})\}\cap Z\Bigr]  \,\,{\displaystyle\bigcap}\,\,\Bigl [ \{F(\bar{x}+\tau u)-Y(\bar{x}+\tau u,\bar{p})\}\cap Z\Bigr ] =\emptyset
.\label{contr'}
\end{equation}
Let us further set:
\[
K_{\tau}:= \{F(\bar{x}+\tau u)-Y(\bar{x}+\tau u)\}\cap Z\text{\qquad
and\qquad}K_{0}:=\{F(\bar{x})-Y(\bar{x})\}\cap Z,
\]
and define, for every $\tau\in\lbrack0,\delta)$ the function $\phi_{\tau
}:\mathbb{R}\mapsto\mathbb{R}\cup\{+\infty\}$ as follows:
\[
\phi_{\tau}(t):=\inf\{\beta\geq0,\;te_{1}+\beta e_{2}\in K_{\tau
}\},\text{\quad}t\in\mathbb{R}\text{.}
\]
Note that $\phi_{\tau}$ is convex with $\phi_{\tau}(0)=0$ and $\phi_{\tau
}(t)>0$ when $t\neq0$. Since $F$ is Lipschitz, we have
\[
\left(  F(\bar{x}+\tau u)-Y(\bar{x}+\tau
u,\bar{p})\right)\,\,  \underset{\tau\to 0}{\longrightarrow}\,\,F(\bar{x})-Y(\bar{x},\bar{p})\qquad\text{(in the Painlev\'e-Kuratowski sense)}
\]
and consequently $\mathrm{epi}(\phi_{\tau})$ converges to $\mathrm{epi}
(\phi_{0})$ as $\tau\rightarrow0$ in the Painlev\'{e}-Kuratowski sense.
Therefore, by Attouch Theorem (see \cite[Theorem 2.1]{CT1998} \emph{e.g.}, or
\cite{Attouch1977}) we get
\begin{equation}
\mathrm{gph\,}\left(  \partial\phi_{\tau
}\right)  \,  \underset{\tau\to 0}{\longrightarrow}\,\mathrm{gph\,}\left(  \partial\phi_{0}\right)  \quad\text{(in the
Painlev\'{e}-Kuratowski sense).}\label{Att}
\end{equation}
If $(-\alpha,\alpha)\in\mathrm{dom\,}\phi_{0}:=\{t\in\mathbb{R}:\mathrm{\,}
\phi_{0}(t)<+\infty\},$ for some $\alpha>0,$ then for all $\alpha_{1}
\in(0,\alpha)$ there exists $\tau_{0}>0$ such that $(-\alpha_{1},\alpha
_{1})\subset\mathrm{dom\,}\phi_{\tau}$, for all $\tau\in(0,\tau_0)$ (thanks to the epiconvergence of $\phi_{\tau}$ to $\phi_{0}$). Let us set
\[
a(\tau):=\langle e_{1},\bar{p}-Y(\bar{x}+\tau u,\bar{p})\rangle,\quad \tau\geq 0.
\]
Recalling that $\sigma_{F(\bar{x}+\tau u)}(\bar{p})=\langle\bar{p},Y(\bar
{x}+\tau u,\bar{p})\rangle$ and $\sigma_{F(\bar{x})}(\bar{p})=\langle\bar
{p},\bar{y}\rangle$ and denoting by $\Lambda$ the Lipschitz constant of
$Y(\cdot,\bar{p})$, we deduce that
\begin{equation}\label{thera}
|a(\tau)|\leq\Lambda\tau|e_{1}|.
\end{equation}
Translating formula~\eqref{contr'} into the coordinates $\{e_1, e_2\}$ given in~\eqref{eq:ada} we readily obtain that for
all $b_1\in\mathrm{B}$ the following inequality holds:
\begin{equation}
\langle-\bar{p},\bar{y}\rangle-\tau\langle\bar{p},w\rangle+\tau\delta
+\langle\bar{p},Y(\bar{x}+\tau u,\bar{p})\rangle\,<\,\phi_{\tau}\left(
a(\tau)+\tau\langle e_{1},w\rangle+\tau\delta b_{1}\right).  \label{contr''}
\end{equation}
Observing that $\langle e_{1},w\rangle =|w|^{2}-\langle\bar{p},w\rangle^{2}>0$ and recalling that $\sigma_{F(\bar{x})}(\bar{p})=\langle\bar{p},\bar{y}\rangle$, we deduce from~\eqref{contr''} that
\[
\phi_{\tau}(a(\tau)+\tau\langle e_{1},w\rangle+\tau\delta b_{1})>\tau\left(
\frac{\sigma_{F(\bar{x}+\tau u)}(\bar{p})-\sigma_{F(\bar{x})}(\bar{p})}{\tau
}-\langle\bar{p},w\rangle+\delta\right).
\]
Therefore we deduce from~\eqref{eq:ada} that for sufficiently small $\tau>0$
\[
\phi_{0}(X_{\tau})\,>\,\frac{\delta}{2}\,\tau\,, \qquad\text{where} \quad X_{\tau}:=a(\tau
)+\tau\langle e_{1},w\rangle+\tau\delta b_{1}\,\neq\,0.
\]
By the mean value
theorem, there exist $Z_{\tau}\in\lbrack0,X_{\tau}]$ and $\xi_{\tau}
\in\partial\phi_{\tau}(Z_{\tau})$ (subdifferential of $\phi_{\tau}$ at~$Z_{\tau}$) such that
\[
\phi_{\tau}(X_{\tau})-\phi_{\tau}(0)=\xi_{\tau}X_{\tau}>\frac{\delta\tau}{2}.
\]
Clearly when $\tau\rightarrow0$ we have $X_{\tau}\rightarrow0$ hence $Z_{\tau
}\rightarrow0$. Moreover since
\[
\frac{1}{\tau}|X_{\tau}|\,\leq\,\Lambda|e_{1}|\,+\,|\langle e_{1},w\rangle
|\,+\,\delta
\]
we have
\[
|\xi_{\tau}|\,\geq\,\,\frac{\delta/2}{\Lambda|e_{1}|+|\langle e_{1},w\rangle|+\delta}\,:=\Delta_{0}\,>\,0.
\]
\smallskip\newline

\textit{Claim}. The set $\{\,|\xi_{\tau}|:\,\,\tau\in(0,\delta)\}$ is bounded.\smallskip
\newline

\noindent Let $(-\alpha_{1},\alpha_{1})\subset\mathrm{dom\,}\phi_{\tau}$ for all
$\tau\in\lbrack0,\delta)$. Since the convex function $\phi_{0}$ is $\kappa
$-Lipschitz on $[-\alpha_{1},\alpha_{1})$ for some constant $\kappa>0,$ we
have $\partial\phi_{0}(t)\subset\lbrack-\kappa,\kappa]$ for all $t\in
(-\alpha_{1},\alpha_{1})$. We now deduce from \eqref{Att} that $\partial
\phi_{\tau}(Z_{\tau})\subset\lbrack-\kappa,\kappa],$ for all $\tau\in
(0,\delta),$ and consequently, we may assume (passing to a subsequence) that
$(\xi_{\tau})_{\tau}$ converges  to some $\xi$ and the sequence $(Z_{\tau}
,\xi_{\tau})\in\mathrm{gph\,}\left(  \partial\phi_{\tau}\right)  $ converges
to $(0,\xi)$ which belongs to $\mathrm{gph\,}\left(  \partial\phi_{0}\right)
$ by \eqref{Att}. We deduce that $|\xi|=\lim|\xi_{\tau}|>\delta_{0}>0$, which
is a contradiction since $\phi_{0}$ is differentiable at $0$ with
$\phi_{0}^{\prime}(0)=0,$ because $\bar{y}=Y[\bar{x},\bar{p})$ is a smooth point
of $\partial F(\bar{x})$. The proof is complete. \hfill$\Box$

\bigskip

\begin{remark} (i).  As we already mentioned at the beginning, differentiability of $F$ at $(x, y)$ is straightforward whenever $y\in \mathrm{int }F(x)$. Theorem~\ref{RadeFint} establishes differentiability also for\emph{a.e.} $y\in \partial F(x)$ (with respect to the Lebesgue measure of $\partial F(x)$, taking into account the bi-Lipschitz homeomorphism with the unit sphere $S^{\ell-1}$). \smallskip\newline
(ii). The proof of Theorem~\ref{RadeFint} uses in a crucial way (namely, in~\eqref{thera}) that the set-valued map $F$ is not only Lipschitz, but isotropically Lipschitz. Although this assumption seems indeed essential, we do not have any counterexample so far.
This could be a topic for further investigations.
\end{remark}

\bigskip

\noindent\textbf{Acknowledgement} Major parts of this work were done during mutual
research visits at the University Bretagne-Occidentale (October 2020, November 2022)
and at the TU Wien (April 2022). In each case the authors express their gratitude to the host institutes
for hospitality. \\
The first author acknowledges support from the Austrian Science Fund (FWF, P-36344-N).



\noindent\rule{4cm}{2pt}

\smallskip

\noindent Aris DANIILIDIS

\medskip

\noindent Institut f\"{u}r Stochastik und Wirtschaftsmathematik, VADOR E105-04
\newline TU Wien, Wiedner Hauptstra{\ss }e 8, A-1040 Wien\medskip\newline(on
leave) DIM-CMM, CNRS IRL 2807 \newline Beauchef 851, FCFM, Universidad de
Chile \medskip\newline\noindent E-mail: \texttt{aris.daniilidis@tuwien.ac.at}
\newline\noindent\texttt{https://www.arisdaniilidis.at/}

\medskip

\noindent Research supported by the grants: \smallskip\newline Austrian
Science Fund (FWF P-36344N) (Austria)\newline CMM FB210005 BASAL funds for
centers of excellence (ANID-Chile)\newline

\vspace{0.5cm}

\noindent Marc QUINCAMPOIX\medskip

\noindent  Univ Brest, CNRS UMR 6205\newline
\noindent Laboratoire de Math\'{e}matiques de Bretagne Atlantique\\
6, avenue Victor Le Gorgen, F-29200 Brest, France

\medskip

\noindent E-mail: \texttt{marc.quincampoix@univ-brest.fr}\newline
\noindent\texttt{http://marc.quincampoix.perso.math.cnrs.fr/}

\end{document}